\documentclass{amsart}
\usepackage{amsmath}
\usepackage{graphicx}
\usepackage{amsfonts}
\usepackage{amssymb}

\newtheorem{theorem}{Theorem}
\theoremstyle{plain}

\newtheorem{lemma}{Lemma}

\newtheorem{remark}{Remark}

\numberwithin{equation}{section}

\catcode`@=11
\def\plainvspace@{\def\vspace##1{\noalign{\vskip##1}}}
\def\plainLet@{\relax\iffalse{\fi\let\\=\cr\iffalse}\fi}
\def\aligned{\,\vcenter\bgroup\plainvspace@\plainLet@\openup\jot\m@th\ialign
  \bgroup \strut\hfil$\displaystyle{##}$&$\displaystyle{{}##}$\hfil\crcr}
\def\endaligned{\crcr\egroup\egroup}
\def\matrix{\,\vcenter\bgroup\plainLet@\plainvspace@
    \normalbaselines
  \m@th\ialign\bgroup\hfil$##$\hfil&&\quad\hfil$##$\hfil\crcr
    \mathstrut\crcr\noalign{\kern-\baselineskip}}
\def\endmatrix{\crcr\mathstrut\crcr\noalign{\kern-\baselineskip}\egroup
                \egroup\,}

\def\cases{\left\{\,\vcenter\bgroup\plainvspace@
     \normalbaselines\openup\jot\m@th
      \plainLet@\ialign\bgroup$\displaystyle{##}$\hfil&\quad$\displaystyle{{}##}$\hfil\crcr
      \mathstrut\crcr\noalign{\kern-\baselineskip}}
\def\endcases{\endmatrix\right.}
\def\fint{{\mathbf -}\!\!\!\!\!\!\int}
\def\sfint{{\mathbf-}\!\!\!\!\!\int}
\begin{document}
\title{Sharp form for improved Moser-Trudinger inequality}
\author  {Yilong Ni and Meijun Zhu}
\address{Department of Mathematics\\
The University of Oklahoma\\
Norman, OK 73019\\ }
\begin{abstract}
We prove that  the improved Moser-Trudinger inequality with
optimal coefficient $\alpha =1/2$ holds for all functions on $S^2$
 with zero moments.
\end{abstract}
\maketitle
\section{Introduction}
The standard Moser-Trudinger-Onofri inequality states
(\cite{mosert}, \cite{On}) that on the standard unit sphere $(S^2,
g_0)$ with the induced metric $g_0$ from $R^3$, for any $u\in
C^1(S^2)$,
$$\frac {1}{4 \pi} \int_{S^2} e^{2u} \le \exp\{
\frac 1{4\pi}\int_{S^2}(|\nabla u|^2+2u)\}, $$ and the equality
holds if and only if $e^{2u} g$ is a metric of constant curvature.
In the study of deforming metrics and prescribing curvatures on
$S^2$, this inequality is often used to control the size and
behavior of a new metric $e^{2u}g_0$  near a concentration point.
With certain ``balance" condition on the metric one would guess
that if the metric concentrates, it should concentrate at more
than one point. Thus it is reasonable to ask whether there is some
small constant $\alpha \in (0, 1)$ and a constant $C_\alpha$ such
that
$$\frac {1}{4 \pi} \int_{S^2} e^{2u} \le C_\alpha \exp\{
\frac 1{4\pi}\int_{S^2}(\alpha |\nabla u|^2+2u)\} $$
holds for those functions satisfying certain balance
condition.

It in fact was first observed by Moser in \cite{moser2} that the
above inequality holds for $\alpha =1/2$ if $u(x)$ is symmetric
with respect with the origin  (that is: $u(x_1, x_2, x_3)=u(-x_1,
-x_2, -x_3).$)  In general, Aubin \cite{aubin} proved that if $$u
\in \Lambda:= \{ f(x) \in H^1(S^2) \ : \ \int_{S^2} e^{2f} x_i=0 \
\mbox{for} \ i=1, 2, 3\},$$ where $\{x_1,x_2,x_3\}$ are the
standard coordinates in ${\mathbf R}^3$, then for any given
constant $\alpha \in (1/2, 1)$, there is a constant $C_\alpha$
such that
\begin{equation}
\frac {1}{4 \pi} \int_{S^2} e^{2u} \le C_\alpha \exp\{ \frac
1{4\pi}\int_{S^2}(\alpha |\nabla u|^2+2u)\}.
\label{add4-18-1}
\end{equation} Later, in their study of
prescribing curvature problem on $S^2$, Chang and Yang
\cite{chang} were able to show that for $\alpha$ close to 1, the
optimal constant for $C_\alpha$ in the above inequality is $1$. On
the other hand, using the standard bubbling sequence, one can see
that Aubin's inequality can not hold if $\alpha< 1/2$. Thus the
immediate question is:

{\it \noindent (I): Is there a constant $C_*$ so that
$$\frac {1}{4 \pi} \int_{S^2} e^{2u} \le C_* \exp\{
\frac 1{4\pi}\int_{S^2}(\frac 12|\nabla u|^2+2u)\}$$ holds for all
$u \in \Lambda$?}

\medskip

If the answer to this question is affirmative, one may continue to
ask

\medskip

{\it \noindent  (II):  what is the optimal constant $C_*$? Is it
$1$?}

\medskip

In this short note, we will give an affirmative answer to the
first question. To answer the second question we need to solve a
partial differential equation. So far we have no clue  how to
solve it. See more details in  Remark at the end of this note.

 Let
$$
O=\{u\in H^1(S^2)\mbox{ : }\int_{S^2}e^{2u}x_i =0,\ i=1,2,3 \mbox{
and }\int_{S^2}u =0\}.
$$
For any $u\in O$, we define functional
$$
I(u)=\frac{1}{2}\fint_{S^2}|\nabla u|^2-\ln\fint_{S^2}e^{2u}
$$
We have the following main result.
\begin{theorem}\label{th1}
There exists constant $C\in\mathbf R$, such that $I(u)>C$ for all
$u\in O$. Moreover $\inf_{u\in O}I(u)$ is attained by some $u\in
O$.
\end{theorem}

It is easy to observe that $\inf_{u\in O}I(u) \le 0.$  But it is
not clear yet whether $\inf_{u\in O}I(u) = 0$ or not. On the other
hand, to our surprise, we are able to show that if  a minimizing
sequence blows up (more details will be given later), then
$\inf_{u\in O}I(u) >0$ (through a dedicated asymptotic analysis).
We thus obtain the existence of the extremal for $\inf_{u\in
O}I(u)$. Similar blow-up analysis is quite standard now (see, for
example,  \cite{DJLW} and \cite {CZ}).

\section{proof of the theorem}
Let $u\in O$.  For any $\epsilon\in(0,1/2)$, we define a perturbed
functional
$$
I_\epsilon(u)=\frac{1}{2(1-\epsilon)}\fint_{S^2}|\nabla u|^2
-\ln\fint_{S^2}e^{2u}.
$$
It follows from Aubin's inequality (\ref{add4-18-1})  that
$$
E_\epsilon=\inf_{u\in O}I_\epsilon(u)>-\infty.
$$
Further, one can show that the infimum is attained by some
$u_\epsilon\in O$. Thus $u_\epsilon$ satisfies the following
Euler-Langrange equation:
$$
-\frac{1}{8\pi(1-\epsilon)}\Delta u_\epsilon=\frac{e^{2u_\epsilon}}
{\int_{S^2}e^{2u_\epsilon}}-\frac1{4\pi}+e^{2u_\epsilon}\sum_{i=1}^3a_ix_i,
$$
where $a_i's$ are Langrange multipliers.

We claim: $a_i=0$ for $i=1, 2, 3.$  The proof of this claim is
along the same line as in \cite{chang}.

Let $v(x)$ be a solution to
$$
\Delta v+he^v=c \ \ \ \ \mbox{on} \ \  S^2.
$$
 Kazdan and Warner  (\cite{kw}) showed that $v$ satisfies
$$\sfint_{S^2} e^v\nabla h\cdot\nabla x_i=(2-c)\sfint_{S^2}
e^vhx_i.$$ Let
$$
c=4(1-\epsilon),\qquad
h=16\pi(1-\epsilon)\left(\frac{1}{\int_{S^2}e^{2u_\epsilon}}+\sum_{i=1}^3
a_ix_i\right),\qquad v=2u_\epsilon,
$$
 we have
$$
\fint_{S^2} e^{2u_\epsilon}\nabla(\sum_{j=1}^3a_jx_j)\cdot\nabla
x_i=(2-4(1-\epsilon)) \fint_{S^2}
e^{2u_\epsilon}\left(\frac{1}{\int_{S^2}e^{2u_\epsilon}}
+\sum_{j=1}^3a_jx_j\right)x_i.
$$
Since $\int_{s^2}e^{2 u_\epsilon} x_i=0$, we know  that
$$
\fint_{S^2} e^{2u_\epsilon}\nabla(\sum_{j=1}^3a_jx_j)\cdot\nabla
x_i=(2-4(1-\epsilon)) \fint_{S^2}
e^{2u_\epsilon}\sum_{j=1}^3a_jx_j x_i.
$$
Multiplying both sides by $a_i$ and summing from $i=1$ to $3$, we
obtain that
$$
\fint_{S^2}
e^{2u_\epsilon}|\nabla(\sum_{j=1}^3a_jx_j)|^2=(2-4(1-\epsilon))
\fint_{S^2} e^{2u_\epsilon}|\sum_{j=1}^3a_jx_j |^2.
$$
It follows from $2-4(1-\epsilon)<0$ that $a_i=0$, $i=1,2,3$.

Therefore $u_\epsilon$ satisfies
\begin{equation}
-\Delta
u_\epsilon=8\pi(1-\epsilon)\left(\frac{e^{2u_\epsilon}}{\int_{S^2}
e^{2u_\epsilon}}-\frac{1}{4\pi}\right). \label{add4-18-2}
\end{equation}

Recall that $u_\epsilon$ is a minimizer of $\inf_{u\in
O}I_\epsilon (u)$, thus, if $\int_{S^2}e^{2u_\epsilon}$ stays
bounded as $\epsilon\to 0$, $\int_{S^2} |\nabla u_\epsilon|^2 \le
C$. Then there exist a subsequence $u_{\epsilon_n}$ converging to
$u_0$ in $H^1(S^2)$. Furthermore, $u_0$ is a minimizer for $I(u)$
and Theorem \ref{th1} follows.

From now on, we  assume that up to a subsequence
$\int_{S^2}e^{2u_\epsilon} \to\infty$ as $\epsilon\to 0$, and will
derive an contradiction. For simplicity, we shall not distinguish
a subsequence $\{\epsilon_i\}$ from the original $\{\epsilon\}$.

Let $v_\epsilon=2u_\epsilon-\ln\int_{S^2}e^{2u_\epsilon}$. Then
$v_\epsilon$ satisfies $\int_{S^2}e^{v_\epsilon}=1$,
$v_\epsilon^a:=\sfint v_\epsilon \to -\infty$ as $\epsilon\to 0$,
and
\begin{equation}
\label{eq}
-\Delta v_\epsilon=16\pi(1-\epsilon)(e^{v_\epsilon}-\frac{1}{4\pi}).
\end{equation}
We first have the $L^q$ estimate for $v_\epsilon$  for any
$q\in[1,2)$:
\begin{equation}
\label{w1q}
||\nabla v_\epsilon||_{L^q(S^2)}<C_q.
\end{equation}
In fact, for any $\varphi\in W^{1,q/(q-1)}(S^2)$ with
$\int_{S^2}\varphi=0$ and
$||\varphi||_{W^{1,q/(q-1)}(S^2)}=1$(thus $\varphi\in
L^{\infty}(S^2)$),
$$
\left|\int_{S^2}\nabla v_\epsilon \nabla\varphi\right|=
\left|\int_{S^2}\Delta v_\epsilon \varphi\right|=
16\pi(1-\epsilon)\left|\int_{S^2}(e^{v_\epsilon}-\frac1{4\pi} )
\varphi\right|\le C.
$$

Since $\int_{S^2}e^{v_\epsilon}=1$, $16\pi(1-\epsilon)e^{v_\epsilon}$
converges in measure to $d\mu$, a positive measure on $S^2$, that is
$\int_{S^2}e^{v_\epsilon}\psi\to \int_{S^2}\psi d\mu$, for any
$\psi\in C^0(S^2)$.  Let
$$
R=\{x\in S^2 \mbox{ : }\exists \psi\in C^0(S^2),0\le\psi\le1,
\psi\equiv 1 \mbox{ around }x,\mbox{ s.t. }\int\psi
d\mu<4\pi \}
$$
be the set of ``regular points", and
$$
S=\{x\in S^2\mbox{ : }\exists x_n\in S^2\mbox{ and
}\{\epsilon_n\}, \mbox{ s.t. } \ \lim_{n\to \infty} x_n=x, \
\mbox{and} \ \lim_{n\to\infty} v_{\epsilon_n}(x_n)\to \infty\},
$$
be the set of ``blow-up points". We need the following lemma of
Brezis and Merle \cite{brezis} to initiate our analysis.
\begin{lemma}
$($Brezis-Merle Lemma$)$ Suppose $u$ satisfies $-\Delta u=f$ on a
bounded domain $\Omega\subset {\mathbf R}^2$ and $u|_{\partial
\Omega}=0$ then for any $\delta\in(0,4\pi)$ there exists a
constant $c(\delta)$ such that
$$
\int_{\Omega}\exp\left(\frac{(4\pi-\delta)|u|}{||f||_{L^1(\Omega)}}\right)
\le c(\delta).
$$
\end{lemma}
Although Brezis-Merle Lemma was originally proved for a bounded domain
$\Omega\subset {\mathbf R}^2$, the same result also holds on any domain
$\Omega\subset S^2$.

Using the above lemma, we derive the following property for
regular points.
\begin{lemma}
\label{bound} For any $x_0\in R$, there exists $r>0$, such that
$v_\epsilon-v_\epsilon^a$ is bounded in $L^\infty(B_r(x_0))$
uniformly in $\epsilon$, where $v_\epsilon^a=\sfint v$ and
$B_r(x_0)$ is the ball centered at $x_0$ with radius $r$.
\end{lemma}
\begin{proof}
By definition, for any $x_0\in R$, there exists $r>0$, such that
$||e^{v_\epsilon}||_{L^1(B_{4r}(x_0))}<4\pi$. Let $v_\epsilon^1$
be the unique solution to the following Dirichlet problem:
$$
\left\{
\begin{array}{rll}
&-\Delta v_\epsilon^1=16\pi(1-\epsilon)e^{v_\epsilon} & \mbox{ in }B_{4r}(x_0)\\
&v_\epsilon^1=0 & \mbox{ on }\partial B_{4r}(x_0),
\end{array}
\right.
$$
and $v_\epsilon^2=v_\epsilon-v_\epsilon^1-v_\epsilon^a$. Applying Brezis-Merle
Lemma to $v_\epsilon^1$, we obtain that $\int_{B_{4r}(x_0)}e^{p|v_\epsilon^1|}
<\infty$, for some $p\in(1,2)$. Since $\Delta v_\epsilon^2=4(1-\epsilon)$,
using $L^p$ interior estimate, we have
\begin{align*}
||v_\epsilon^2||_{L^\infty(B_{2r}(x_0))}&\le C||v_\epsilon^2||_{W^{2,p}
(B_{2r}(x_0))}\le C||v_\epsilon^2||_{L^p(B_{4r}(x_0))}\\
&\le C(||v_\epsilon-v_\epsilon^a||_{L^p(B_{4r}(x_0))}
+||v_\epsilon^1||_{L^p(B_{4r}(x_0))})\\
&\le C(||\nabla v_\epsilon||_{L^p(S^2)}
+||v_\epsilon^1||_{L^p(B_{4r}(x_0))})\\
&\le C
\end{align*}
where we also use (\ref{w1q}) and Poincar\'e inequality. Therefore
$$
\int_{B_{2r}(x_0)}e^{pv_\epsilon}=\int_{B_{2r}(x_0)}e^{pv_\epsilon^a}\cdot
e^{pv_\epsilon^1}\cdot e^{pv_\epsilon^2}\le C.
$$
It follows easily from interior $L^p$ estimate (for $p>1$)  that
$||v_\epsilon^1||_{L^{\infty}(B_r(x_0))}\le C$. Hence
$||v_\epsilon-v_\epsilon^a||_{L^{\infty}(B_r(x_0))}\le C$.
\end{proof}

Using the above lemma, we immediately get that $S\subset S^2\setminus R$.
Therefore
$$
\#S\le \#( S^2\setminus R)\le \frac{\int_{S^2} d\mu}{4\pi}\le 4,
$$
where $\#S$ is the cardinality of $S$. Choose $r$ small, so that
for $x\in S$, $B_r(x)$ are disjoint. For any $x\in S$, by
definition, for $\epsilon$ sufficiently small, $v_\epsilon$ has a
local maximum $x_\epsilon\in B_r(x)$ with $v_\epsilon (x_\epsilon)
\to \infty$, and up to a subsequence $x_\epsilon\to x$ as
$\epsilon\to 0$. Choose a normal coordinate system around $x$ and
define
$$
\varphi_\epsilon(x)=v_\epsilon(\tau_\epsilon^{-1}x+x_\epsilon)
-\lambda_\epsilon,
$$
where $\lambda_\epsilon=v_\epsilon(x_\epsilon)$, $\tau_\epsilon
=e^{\lambda_\epsilon/2}$ and we use
$\tau_\epsilon^{-1}x+x_\epsilon$ to represent
$\exp_{x_\epsilon}(\tau_\epsilon^{-1}x)$. For fixed $R>0$, when
$\epsilon$ is sufficiently small, $\varphi_\epsilon$ satisfies
$$
-\Delta \varphi_\epsilon=16\pi(1-\epsilon)(e^{\varphi_\epsilon}
-\frac{1}{4\pi\tau_\epsilon^2}) \ \mbox{ in } \
B_{2R}(0)\subset{\mathbf R}^2.
$$
\begin{lemma}
For a fixed $R>0$, $\varphi_\epsilon$ is bounded in $B_R(0)$
uniformly in $\epsilon$.
\end{lemma}
\begin{proof}
Let $\varphi_\epsilon^{(1)}$ be the unique solution to
$$
\cases &-\Delta
\varphi_\epsilon=16\pi(1-\epsilon)(e^{\varphi_\epsilon}
-\frac{1}{4\pi\tau_\epsilon^2}) \ \mbox{ in } \ B_{2R}(0)\subset{\mathbf R}^2.\\
&\varphi^{(1)}_\epsilon|_{\partial B_{2R}(0)}=0.
\endcases
$$
Since $x_\epsilon$ is a local maximum point of $v_\epsilon(x)$, we
have $\varphi_\epsilon\le \varphi_\epsilon(0)=0$ and $e^{\varphi_\epsilon}\le 1$. It
follows that $||\varphi_\epsilon^{(1)}||_{L^\infty}\le C<+\infty$.
Let $\varphi_\epsilon^{(2)}=\varphi_\epsilon-\varphi_\epsilon^{(1)}$. Then
$\varphi_\epsilon^{(2)}\le -\varphi_\epsilon^{(1)}\le C$. Since
$2C-\varphi_\epsilon^{(2)}\ge C$ is harmonic, positive and
$$
2C-\varphi_\epsilon^{(2)}(0)=2C-\varphi_\epsilon(0)+\varphi^{(1)}_\epsilon(0)
 \le 3C,
$$
Harnack's inequality implies that
$||2C-\varphi_\epsilon^{(2)}||_{L^\infty}\le \tilde{C}$ in $B(R)$. Hence
$||\varphi_\epsilon ||_{L^\infty(R)}\le C+\tilde{C}$.
\end{proof}

Since $\varphi_\epsilon$ is uniformly bounded in $B_R(0)$,
elliptic estimates yield that, up to a subsequence, $\varphi_\epsilon \to
\varphi_0$ in $C^{2,\alpha}(B(R/2))$ for some $\alpha\in(0,1)$ and
$\varphi_0$ satisfies
$$
\cases &-\Delta \varphi_0=16\pi e^{\varphi_0} \ \mbox{
in }{\mathbf R}^2\\
&\varphi_0(0)=0. \endcases
$$
Furthermore, $\varphi_\epsilon(x)\le \varphi_\epsilon(0)=0$ and
$$
\int_{{\mathbf R}^2}e^{\varphi_0}\le \overline{\lim_{R\to\infty}}\,
\overline{\lim_{n\to\infty}} \int_{B(R)}e^{\varphi_{\epsilon_n}}\le1.
$$
The uniqueness theorem in \cite{chen} implies that
$$
\varphi_0(x)=2\ln\frac{1}{1+2\pi |x|^2}.
$$
Therefore, as $\epsilon\to 0$,
$$
1=\int_{S^2}e^{v_\epsilon}\ge \sum_{x\in
S}\int_{B_r(x)}e^{v_\epsilon} \ge \sum_{x\in
S}\int_{B_R(0)}e^{\varphi_\epsilon} \to\sum_{x\in S}\frac{\pi
R^2}{1+2\pi R^2}.
$$
It follows that $\#S\le 2$ and $R=S^2\setminus S$. Since $S$ is
not empty and $v_\epsilon\in O$, we obtain that $\#S=2$ and
$S=\{{\bf a,-a}\}$ for some ${\bf a}\in S^2$. Without loss of
generality, we may assume that $S=\{n,s\}$, where $n$ and $s$
stand for the north pole and the south pole respectively.

For any compact domain $K\subset\subset S^2\setminus S(=R)$, we know from
Lemma \ref{bound} that $v_\epsilon-v_\epsilon^a$ is bounded in
$L^\infty(K)$ uniformly in $\epsilon$. Since $v_\epsilon^2\to -\infty$, it
follows from the standard elliptic estimate that
$$
v_\epsilon-v_\epsilon^a\to G(x) \mbox{ in }C^{1,\alpha}(K),
$$
where $G(x)$ satisfies
$$
\cases
& -\Delta G+4=8\pi(\delta_n+\delta_s)\mbox{ on }S^2, \\
&\int_{S^2}G=0,
\endcases
$$
and $\delta_n$, $\delta_s$ are delta functions centered at the north pole and
the south pole, respectively. It can be easily seen that
$$
G(x)=-4\ln\sin\theta-4(1-\ln 2),
$$
where $\theta$ is the angle between $x$ and $x_3$.

Let $x_{\epsilon 1}$ and $x_{\epsilon2}$ be the local maximum near the north
pole and the south pole respectively. We use the following notations:
\begin{align*}
&\lambda_{\epsilon1}:=v_\epsilon(x_{\epsilon1}), \quad
\lambda_{\epsilon2}:=v_\epsilon(x_{\epsilon2}), \quad
\tau_{\epsilon1}:=e^{\frac{\lambda_{\epsilon1}}2},\quad
\tau_{\epsilon2}:=e^{\frac{\lambda_{\epsilon2}}2},\\
&r_{\epsilon1}:=\frac R{\tau_{\epsilon1}},\quad
r_{\epsilon2}:=\frac R{\tau_{\epsilon2}}, \mbox{ for fixed }R>0,\\
&B_1:=B_{r_{\epsilon1}}(n), \quad B_2:=B_{r_{\epsilon2}}(s),\quad
\Omega:=S^2\setminus(B_1\cup B_2).
\end{align*}
It follows from Theorem 0.2 in \cite{li} that there exists some constant
$C'$, such that
$$
\left|v_\epsilon(x)-\ln\frac{e^{\lambda_{\epsilon 1}}}{(1+2\pi(1-\epsilon)
e^{\lambda_{\epsilon 1}}\mbox{dist}(n,x)^2)^2}\right|\le C',
\mbox{ for }x\in B_{\pi/4}(n)
$$
and
$$
\left|v_\epsilon(x)-\ln\frac{e^{\lambda_{\epsilon 2}}}{(1+2\pi(1-\epsilon)
e^{\lambda_{\epsilon 2}}\mbox{dist}(s,x)^2)^2}\right|\le C',
\mbox{ for }x\in B_{\pi/4}(s).
$$
Since $v_\epsilon-v_\epsilon^\alpha\to G$ in
$C^{1,\alpha}(S^2\setminus[B_{\pi/4}(n) \cup B_{\pi/4}(s)])$, we
obtain that there exists some constant $C$, such that
\begin{equation}\label{yyli}
|\lambda_{\epsilon1}-\lambda_{\epsilon2}|\le C,
\end{equation}
for all $\epsilon$. Without loss of generality, we may assume that
$\lambda_{\epsilon1}\ge\lambda_{\epsilon2}$.
\begin{lemma}\label{lde}
For all $x\in\Omega$, $v_\epsilon(x)\ge G(x)+D_\epsilon+o_\epsilon(1)$,
where
\begin{equation}
\label{depsilon}
D_\epsilon=-\lambda_{\epsilon1}+2\ln\frac{R^2}{1+2\pi R^2}+4(1-\ln2)
\end{equation}
and $o_\epsilon(1)$ stands for some function that goes to $0$ as
$\epsilon\to 0$.
\end{lemma}

\begin{proof}
In $B_1$, we have
\begin{align*}
G(x)&=-4\ln r-4(1-\ln2)+o_\epsilon(1),\\
v_\epsilon(x)&=\lambda_{\epsilon1}+2\ln\frac{1}{1+2\pi
\tau^2_{\epsilon1}r^{2}}+o_\epsilon(1).
\end{align*}
Hence,
\begin{align*}
(v_\epsilon-G)|_{\partial B_1}&=\lambda_{\epsilon1}+2\ln\frac{1}
{1+2\pi R^2}+4\ln\frac{R}{\tau_{\epsilon1}}+4(1-\ln2)+o_\epsilon(1)\\
=&-\lambda_{\epsilon1}+2\ln\frac{R^2}{1+2\pi R^2}+4(1-\ln2)+o_\epsilon(1)\\
=&D_\epsilon+o_\epsilon(1).
\end{align*}
Similarly
\begin{align*}
(v_\epsilon-G)|_{\partial B_2}&=\lambda_{\epsilon2}+2\ln\frac{1}
{1+2\pi R^2}+4\ln\frac{R}{\tau_{\epsilon2}}+4(1-\ln2)+o_\epsilon(1)\\
=&-\lambda_{\epsilon2}+2\ln\frac{R^2}{1+2\pi R^2}+4(1-\ln2)+o_\epsilon(1)\\
=&\lambda_{\epsilon1}-\lambda_{\epsilon2}+D_\epsilon+o_\epsilon(1)\ge
D_\epsilon+o_\epsilon(1).
\end{align*}
Since $\Delta(v_\epsilon-G)\le 0$ in $\Omega$ and $(v_\epsilon-G)
|_{\partial\Omega}\ge D_\epsilon+o_\epsilon(1)$, Lemma \ref{lde} follows from
the maximum principle.
\end{proof}
We are now ready to estimate $E_\epsilon=I_\epsilon(u_\epsilon)$.
$$
\int_{S^2}|\nabla v_\epsilon|^2 =\left(\int_{B_1}
+\int_{B_2}+\int_{\Omega}\right)|\nabla v_\epsilon|^2
:=I_1+I_2+I_3.
$$
Since the behavior of $v_\epsilon$ near the north pole can be described by the
behavior of $\phi_\epsilon$ in $B_{R}(0)\subset {\mathbf R}^2$ and
$\phi_\epsilon\to\phi_0$ in $C^{1,\alpha}$, we obtain that
\begin{align*}
I_1&=\int_{B_1}|\nabla v_\epsilon|^2=\int_{B_R(0)}
\left|\nabla\left(2\ln\frac{1}{1+2\pi |x|^2}\right)\right|^2dx+o_\epsilon(1)\\
&=128\pi^3\int_0^R\frac{r^2rdr}{(1+2\pi r^2)^2}
=16\pi(\ln(1+2\pi R^2)-1)+o_\epsilon(1)+o_R(1),
\end{align*}
where $o_R(1)\to 0$ as $R\to\infty$. Similarly
$I_2=16\pi(\ln(1+2\pi R^2)-1)+o_\epsilon(1)+o_R(1)$.

For $I_3$ we have
$$
I_3=\int_\Omega |\nabla v_\epsilon|^2=-\int_\Omega
v_\epsilon\Delta v_\epsilon
+\int_{\partial\Omega}v_\epsilon\frac{\partial
v_\epsilon}{\partial\mathbf n}.
$$
Using Green's formula, we obtain that
\begin{align*}
-\int_\Omega v_\epsilon\Delta v_\epsilon=&\int_\Omega((-\Delta v_\epsilon
+4(1-\epsilon))v_\epsilon-4(1-\epsilon)v_\epsilon\\
\ge&\int_\Omega(-\Delta v_\epsilon+4(1-\epsilon))(G+D_\epsilon)
-4(1-\epsilon)v_\epsilon+o_\epsilon(1)\\
=&-\int_\Omega\Delta v_\epsilon(G+D_\epsilon)+4(1-\epsilon)\int_\Omega
(G+D_\epsilon-v_\epsilon)+o_\epsilon(1)\\
=&-\int_\Omega v_\epsilon\Delta G+\int_{\partial \Omega}(v_\epsilon
\frac{\partial G}{\partial\mathbf n}-(G+D_\epsilon)
\frac{\partial v_\epsilon}{\partial\mathbf n})\\
&\hskip 2cm +4(1-\epsilon)\int_\Omega
(G+D_\epsilon-v_\epsilon)+o_\epsilon(1).
\end{align*}
It follows from $\Delta G=4$ on $\Omega$, that
$$
I_3\ge -4(2-\epsilon)\int_\Omega v_\epsilon+4(1-\epsilon)\int_\Omega(G+D_\epsilon)
+\int_{\partial\Omega}(v_\epsilon\frac{\partial G}{\partial\mathbf n}
+(v_\epsilon-G-D_\epsilon)\frac{\partial v_\epsilon}{\partial\mathbf n}).
$$
We now estimate the terms in the right hand side of the above inequality. First it
is easy to see that
$$
\int_\Omega G=\left(\int_{S^2}-\int_{B_1\cup B_2}\right)G
=0-\int_{B_1\cup B_2}G=o_\epsilon(1)
$$
and
$$
\int_\Omega v_\epsilon=\int_{S^2} v_\epsilon-\int_{B_1\cup B_2} v_\epsilon
=4\pi v_\epsilon^a+o_\epsilon(1).
$$
Since $\lambda_{\epsilon i}\cdot\mbox{vol}(B_i)=o_\epsilon(1)$, $i=1,2$,
it follows from (\ref{yyli}) and (\ref{depsilon}) that $D_\epsilon\cdot
 \mbox{vol}(B_i)=o_\epsilon(1)$, $i=1,2$. Therefore
$$
\int_\Omega D_\epsilon=\int_{S^2} D_\epsilon-\int_{B_1\cup B_2} D_\epsilon
=4\pi D_\epsilon+o_\epsilon(1).
$$
On $\partial B_i$ for $i=1,2$
\begin{align*}
G(x)& =4\ln\frac1{{r_{\epsilon i}}}-4(1-\ln2)+o_\epsilon(1),\\
\frac{\partial G}{\partial{\mathbf n}}&=-4\frac{1}
{r_{\epsilon i}}+o_\epsilon(1),\\
v_\epsilon(x)&=\lambda_{\epsilon i}+2\ln\frac{1}{1+2\pi R^{2}}+o_\epsilon(1),\\
\frac{\partial v_\epsilon}{\partial{\mathbf n}}&=(-\frac{8\pi R
}{1+2\pi R^{2}}+o_\epsilon(1)) \tau_{\epsilon i},
\end{align*}
It follows that
\begin{align*}
\int_{\partial\Omega}v_\epsilon\frac{\partial G}{\partial\mathbf n}
=&-\left(\int_{\partial B_1}+\int_{\partial B_2}
\right)v_\epsilon\frac{\partial G}{\partial\mathbf n}\\
=&2\pi r_{\epsilon1}(\frac{4}{r_{\epsilon1}}+o_\epsilon(1))
(\lambda_{\epsilon1}+2\ln\frac{1}{1+2\pi R^2}+o_{\epsilon}(1))\\
&+2\pi r_{\epsilon2}(\frac{4}{r_{\epsilon2}}+o_\epsilon(1))
(\lambda_{\epsilon2}+2\ln\frac{1}{1+2\pi R^2}+o_{\epsilon}(1))\\
=&8\pi(\lambda_{\epsilon1}+\lambda_{\epsilon2}) -32 \pi \ln(1+2\pi
R^2)+o_\epsilon(1)
\end{align*}
and
\begin{align*}
\int_{\partial\Omega}(v_\epsilon-G-D_\epsilon)
\frac{\partial v_\epsilon}{\partial\mathbf n}&=
-\left(\int_{\partial B_1}+\int_{\partial B_2}\right)
(v_\epsilon-G-D_\epsilon)\frac{\partial v_\epsilon}{\partial\mathbf n}\\
&=8\pi(\lambda_{\epsilon1}-\lambda_{\epsilon2})+o_\epsilon(1)+o_R(1).
\end{align*}
Using estimates of $I_1$, $I_2$ and $I_3$, we finally obtain that
\begin{align*}
\int_{S^2}|\nabla v_\epsilon|^2=&I_1+I_2+I_3\ge 32\pi(\ln(1+2\pi R^2)-1)
-4(2-\epsilon)4\pi v_\epsilon^a\\
&\quad  +8\pi(\lambda_{\epsilon1}+\lambda_{\epsilon2})-32\pi\ln(1+2\pi R^2)
+8\pi(\lambda_{\epsilon1}-\lambda_{\epsilon2})\\
&\quad  +16\pi(1-\epsilon)
D_\epsilon+o_\epsilon(1)+o_R(1)\\
=&-32\pi-16\pi(2-\epsilon)v_\epsilon^a+8\pi(\lambda_{\epsilon1}
+\lambda_{\epsilon2})+8\pi(\lambda_{\epsilon1}-\lambda_{\epsilon2})\\
&\quad  +16\pi(1-\epsilon)\left(-\lambda_{\epsilon1}+2\ln\frac{R^2}{1+2\pi R^2}
  +4(1-\ln 2)\right)\\
&\quad  +o_\epsilon(1)+o_R(1)\\
\ge&32\pi-16\pi(2-\epsilon)v_\epsilon^a -32\pi\ln(8\pi)+o_\epsilon(1)+o_R(1).
\end{align*}
It follows that
\begin{align*}
I_\epsilon(u_\epsilon)&=\frac{1}{2(1-\epsilon)}\fint_{S^2}
\frac{|\nabla v_\epsilon|^2}{4}-\ln\fint _{S^2} e^{2u_\epsilon}\\
&\ge 1-\frac{2-\epsilon}{2(1-\epsilon)}v_\epsilon^a
 -\ln(8\pi) +\ln(4\pi)+v_\epsilon^a+ o_\epsilon(1)+o_R(1)\\
&\ge 1-\ln 2 +o_\epsilon(1)+o_R(1).
\end{align*}
Hence $I(u)\ge 1-\ln 2$ for all $u\in O$. This contradicts with the
fact that $I(0)=0$. Therefore $\int_{S^2}e^{2u_\epsilon}$ stays bounded
as $\epsilon\to 0$ and Theorem \ref{th1} follows.

\medskip

\begin{remark}
It is clear from the  proof our main theorem that the extremal
function $u \in O$ satisfies
\begin{equation}
-\Delta u =8\pi\left(\frac{e^{2u }}{\int_{S^2} e^{2u
}}-\frac{1}{4\pi}\right) \ \ \ \ \ \mbox{on} \ \ \ \ S^2.
\label{add4-25-1}
\end{equation}
To answer question (II) in the introduction, one needs to answer
whether $u=0$ is the only solution to (\ref{add4-25-1}). We do not
know the answer yet.
\end{remark}


\medskip
\begin{thebibliography}{9}

\bibitem{aubin} T. Aubin, Meilleures constantes dans le th¨¦or¨¨me d'inclusion de
Sobolev et un th¨¦or¨¨me de Fredholm non lin¨¦aire pour la
transformation conforme de la courbure scalaire,q (French) J.
Funct. Anal. 32 (1979), no. 2, 148--174.

\bibitem{brezis} H. Brezis, F. Merle, Uniform estimates and blow up behavior
for solutions of $-\Delta u=V(x)e^u$ in two dimensions, Comm. Partial Diff.
Equation 16(1991), 1223-1253.

\bibitem{chang} S. Y. A. Chang and P. Yang, Prescribing Gaussian curvature on
$S\sp 2$, Acta Math. 159 (1987), no. 3-4, 215-259.

\bibitem{chen} W. Chen and C. Li, Classification of solutions of some nonlinear
elliptic equations, Duke Math. J. 63(1991). 615-622.

\bibitem{CZ}X.X. Chen and M. Zhu,  Liouville energy on a topological two
sphere, preprint.

\bibitem{DJLW}Ding, W., Jost, J., Li, J. and  Wang, G.   An analysis of
the two-vortex case in the Chern-Simons Higgs model. Calc. Var.
Partial Differential Equations 7 (1998), no. 1, 87--97.


\bibitem{kw} J. Kazdan and F. Warner, Curvature functions for open
$2$-manifolds. Ann. of Math. (2) 99 (1974), 203-219.

\bibitem{zhu}J. Li  and M. Zhu, Sharp local embedding inequalities,
Comm. Pure Appl. Math. 59(2006), 122-144.

\bibitem{li} Y. Y. Li, Harnack type inequality: the method of moving planes,
Comm. Math. Phys. 200(1999), no. 2, 421--444.

\bibitem{mosert} J. Moser, A sharp form of an inequality by N. Trudinger,
Indiana Univ. Math. J. 11(1971), 1077-1092.

\bibitem{moser2} J. Moser, On a nonlinear problem in differential
geometry. Dynamical systems (Proc. Sympos., Univ. Bahia, Salvador,
1971), pp. 273--280.

 \bibitem{On}E. Onofri, On the positivity of the effective action
 in a theory of random surfaces. Comm. Math. Phys. 86 (1982), no. 3,
 321--326.

\bibitem{trudinger} N. Trudinger, On imbeddings into Orlicz spaces and some
applications, J. Math. Mech. 17(1967), 473-483.

\end{thebibliography}
\end{document}